\documentclass[12pt]{article}
\usepackage{amsmath,amssymb,amsthm}
\numberwithin{equation}{section}
\usepackage{units}
\usepackage{color}
\usepackage[T1]{fontenc}
\usepackage[utf8]{inputenc}
\usepackage{authblk}
\usepackage{bm}
\usepackage{graphicx}

\usepackage{datetime}

\usepackage[colorlinks=true,
linkcolor=webgreen,
filecolor=webbrown,
citecolor=webgreen]{hyperref}

\definecolor{webgreen}{rgb}{0,.5,0}
\definecolor{webbrown}{rgb}{.6,0,0}

\usepackage[toc,page]{appendix}

\newcommand{\Z}{{\mathbb Z}}

\newtheorem{thm}{Theorem}
\newtheorem{theorem}[thm]{Theorem}
\newtheorem{lemma}{Lemma}

\newtheorem{cor}[thm]{Corollary}

\title{Linear Properties of Generalized $n$-step Fibonacci Numbers}

\author[]{Kunle Adegoke \\\href{mailto:adegoke00@gmail.com}{\tt adegoke00@gmail.com}}

\affil{Department of Physics and Engineering Physics, \mbox{Obafemi Awolowo University}, 220005 Ile-Ife, Nigeria}

\begin{document}

\date{}

\maketitle

\begin{abstract}
\noindent
We present numerous interesting, mostly new, results involving the $n$-step Fibonacci numbers and $n$-step Lucas numbers and a generalization. Properties considered include recurrence relations, summation identities, including binomial and double binomial summation identities, partial sums and ordinary generating functions. Explicit examples are given for small $n$ values.
\end{abstract}
\section{Introduction}
For $n\ge 2$, the $n$-step Fibonacci numbers, $U_r$ ($r\ge n$), satisfy the linear recurrence relation~\cite{miles60,noe05,howard11}
\begin{equation}\label{eq.qnse2iy}
U_r  = U_{r - 1}  + U_{r - 2}  + U_{r - 3}  +  \cdots  + U_{r - n}  = \sum_{i = 1}^n {U_{r - i} }\,,
\end{equation}
with $n$ initial terms
\begin{equation}\label{eq.c8583mo}
U_k  = 0,\quad-n+2\le k\le 0,\quad U_{-n+1}  =  1\,.
\end{equation}
Well-known members of this number family include the Fibonacci numbers $F_r$ ($n=2$, $U=F$), the Tribonacci numbers $T_r$ ($n=3$, $U=T$), the Tetranacci numbers $M_r$ ($n=4$, $U=M$). The reader is referred to Table~\ref{tab.afmsmgs} for notation and nomenclature.

\medskip

By writing $U_{r-1}  = U_{r - 2}  + U_{r - 3}  + U_{r - 4}  +  \cdots  + U_{r - n -1} $ and substracting this from relation~\eqref{eq.qnse2iy}, we see that the $n$-step Fibonacci numbers also obey the following recurrence relation:
\begin{equation}\label{eq.iiy0jzg}
U_r=2U_{r-1}-U_{r-n-1}\,.
\end{equation} 

\medskip

Extension of the definition of $n$-step Fibonacci numbers to negative subscripts $r<-n+2$ is provided by writing the recurrence relation~\eqref{eq.iiy0jzg} as
\begin{equation}\label{eq.x1xp283}
U_{ - r}  = 2U_{ - r + n}  - U_{ - r + n + 1}\,.
\end{equation}
From \eqref{eq.qnse2iy}, \eqref{eq.c8583mo}, \eqref{eq.iiy0jzg} and \eqref{eq.x1xp283}, we have the following special values:
\begin{equation}\label{eq.fz5izyy}
U_1  = 1,\quad U_k  = \sum_{j = 1}^{k - 1} {U_j } ,\quad 2 \le k \le n - 1,\quad U_{-1}=\delta_{n,2}, \quad U_{ - n}  =  - 1,\quad U_{ - n - 1}  = 2\delta_{n,2}\,,
\end{equation}
where $\delta_{i,j}$ is Kronecka delta, equals $1$ when $i=j$ and equals $0$ otherwise.

\medskip

We also have
\begin{equation}
U_n  = 2^{n - 2} ,\quad U_{n + 1}  = 2^{n - 1} ,\quad U_{n + 2}  = 2^n  - 1\,,
\end{equation}
and, in fact,
\begin{equation}\label{eq.d011uml}
U_{n + k}  = 2^{n + k - 2}  - \sum_{j = 1}^k {2^{j - 1} U_{k - j} },\quad k\in\Z\,.
\end{equation}
We remark that identity \eqref{eq.d011uml} is equivalent to Theorem 3.1 of Howard and Cooper \cite{howard11} without a restriction on $k$. Note that identity \eqref{eq.d011uml} is a special case of identity \eqref{eq.noe5wpc}.

\medskip

Like the $n$-step Fibonacci numbers, the $n$-step Lucas numbers \cite{noe05} obey an $n$th order recurrence relation
\begin{equation}\label{eq.i5d57xw}
V_r  = V_{r - 1}  + V_{r - 2}  + V_{r - 3}  +  \cdots  + V_{r - n}  = \sum_{i = 1}^n {V_{r - i} }\,,
\end{equation}
but with the initial terms
\begin{equation}\label{eq.p5wdyuq}
V_k  = -1,\quad-n+1\le k\le -1,\quad V_0  =  n\,.
\end{equation}
The most well-known members of the $n$-step Lucas numbers are the Lucas numbers ($n=2$), $(L_r)_{r\in\Z}$, and the Tribonacci-Lucas numbers ($n=3$), $(K_r)_{r\in\Z}$.

\medskip

The $n$-step Lucas numbers also obey the three-term recurrence relation
\begin{equation}\label{eq.t3gnmwy}
V_r=2V_{r-1}-V_{r-n-1}\,.
\end{equation}
Extension of the definition of $n$-step Lucas numbers to integers $r<-n+1$ is provided through
\begin{equation}\label{eq.apavzx3}
V_{ - r}  = 2V_{ - r + n}  - V_{ - r + n + 1}\,.
\end{equation}
Noe and Post~\cite{noe05} noted that the $n$-step Fibonacci numbers and the $n$-step Lucas numbers are connected through the identity
\begin{equation}\label{eq.j032ewx}
V_r  = U_r  + 2U_{r - 1}  +  \cdots  + (n - 1)U_{r - n + 2}  + nU_{r - n + 1}  = \sum_{j = 1}^{n} {jU_{r - j + 1} }\,.
\end{equation}
From identities~\eqref{eq.qnse2iy}, \eqref{eq.iiy0jzg} and \eqref{eq.j032ewx}, we can derive the following four-term relation 
\begin{equation}
V_r  = V_{r - 1}  - (n + 1)U_{r - n}  + 2U_r\,,
\end{equation}
which can also be written in the alternative form
\begin{equation}\label{eq.q6egv82}
V_r  = V_{r - 1}  - nU_{r - n}  + U_{r + 1}
\end{equation}
or
\begin{equation}
V_r  = V_{r - 1}  - 2nU_r  + (n + 1)U_{r + 1}\,.
\end{equation}
From \eqref{eq.i5d57xw}, \eqref{eq.p5wdyuq}, \eqref{eq.t3gnmwy}, \eqref{eq.apavzx3} and \eqref{eq.q6egv82}, we also have the following special values for the $n$-step lucas numbers:
\begin{equation}
V_1  = 1,\quad V_{ - n}  = 2n - 1,\quad V_{ - n - 1}  =  - n - 2,\quad V_n  = 2^n  - 1\,.
\end{equation}
The first few sequences of the $n$-step Fibonacci numbers and the $n$-step Lucas numbers are presented in Table~\ref{tab.zcsobmk}.

\medskip

The generalized $n$-step Fibonacci numbers, $W_r$, satisfy the same recurrence equation given in~\eqref{eq.qnse2iy} but with arbitrary initial values. Thus, 
\begin{equation}\label{eq.ah0hczd}
W_r  = W_{r - 1}  + W_{r - 2}  + W_{r - 3}  +  \cdots  + W_{r - n}  = \sum_{i = 1}^n {W_{r - i} }\,,
\end{equation}
for $r\ge n$ but $W_0$, $W_1$, $\ldots$, $W_{n-1}$ are arbitrary. Analogous to~\eqref{eq.iiy0jzg} and \eqref{eq.x1xp283}, we have
\begin{equation}\label{eq.b5q12ro}
W_r=2W_{r-1}-W_{r-n-1}
\end{equation} 
and
\begin{equation}\label{eq.op5kp4d}
W_{ - r}  = 2W_{ - (r - n)}  - W_{ - (r - n - 1)}\,.
\end{equation}
\begin{table}[h!]\label{tab.afmsmgs}
\begin{tabular}{llllllll}
$n$ & Name & Symbol &  &  & $n$ & Name & Symbol \\ 
\cline{1-3}\cline{6-8}
$2$ & Fibonacci & \multicolumn{1}{c}{$F$} &  &  & $6$ & Sextanacci & \multicolumn{1}{c}{$S$} \\ 
 & Fibonacci-Lucas & \multicolumn{1}{c}{$L$} &  &  &  & Sextanacci-Lucas & \multicolumn{1}{c}{} \\ 
 & Generalized Fibonacci & \multicolumn{1}{c}{$\mathcal{F}$} &  &  &  & Generalized Sextanacci & \multicolumn{1}{c}{$\mathcal{S}$} \\ 
\cline{1-3}\cline{6-8}
$3$ & Tribonacci & \multicolumn{1}{c}{$T$} &  &  & $7$ & Heptanacci & \multicolumn{1}{c}{$H$} \\ 
 & Tribonacci-Lucas & \multicolumn{1}{c}{$K$} &  &  &  & Heptanacci-Lucas & \multicolumn{1}{c}{} \\ 
 & Generalized Tribonacci & \multicolumn{1}{c}{$\mathcal{T}$} &  &  &  & Generalized Heptanacci & \multicolumn{1}{c}{$\mathcal{H}$} \\ 
\cline{1-3}\cline{6-8}
$4$ & Tetranacci & \multicolumn{1}{c}{$M$} &  &  & $8$ & Octanacci & \multicolumn{1}{c}{$O$} \\ 
 & Tetranacci-Lucas & \multicolumn{1}{c}{$R$} &  &  &  & Octanacci-Lucas & \multicolumn{1}{c}{} \\ 
 & Generalized Tetranacci & \multicolumn{1}{c}{$\mathcal{M}$} &  &  &  & Generalized Octanacci & \multicolumn{1}{c}{$\mathcal{O}$} \\ 
\cline{1-3}\cline{6-8}
$5$ & Pentanacci & \multicolumn{1}{c}{$P$} &  &  & $9$ & Nanonacci & \multicolumn{1}{c}{$N$} \\ 
 & Pentanacci-Lucas & \multicolumn{1}{c}{$Q$} &  &  &  & Nanonacci-Lucas &  \\ 
 & Generalized Pentanacci & \multicolumn{1}{c}{$\mathcal{P}$} &  &  &  & Generalized Nanonacci & \multicolumn{1}{c}{$\mathcal{N}$} \\ 
\end{tabular}
\caption{Notation and nomenclature for the first few members of the $n$-step Fibonacci numbers, $n$-step Lucas numbers and the generalized $n$-step Fibonacci numbers.}
\end{table}
\begin{table}[h!]\label{tab.zcsobmk}
\begin{tabular}{llllllllllllllllll}
$n$ & Name & $r$ & \multicolumn{1}{c}{$-4$} & \multicolumn{1}{c}{$-3$} & \multicolumn{1}{c}{$-2$} & \multicolumn{1}{c}{$-1$} & \multicolumn{1}{c}{$0$} & \multicolumn{1}{c}{$1$} & \multicolumn{1}{c}{$2$} & \multicolumn{1}{c}{$3$} & \multicolumn{1}{c}{$4$} & \multicolumn{1}{c}{$5$} & \multicolumn{1}{c}{$6$} & \multicolumn{1}{c}{$7$} & \multicolumn{1}{c}{$8$} & \multicolumn{1}{c}{$9$} & \multicolumn{1}{c}{$10$} \\ 
\hline
$2$ & Fibonacci & $F_r$ & \multicolumn{1}{c}{$ -3$} & \multicolumn{1}{c}{$2$} & \multicolumn{1}{c}{$ -1$} & \multicolumn{1}{c}{$1$} & \multicolumn{1}{c}{$0$} & \multicolumn{1}{c}{$1$} & \multicolumn{1}{c}{$1$} & \multicolumn{1}{c}{$2$} & \multicolumn{1}{c}{$3$} & \multicolumn{1}{c}{$5$} & \multicolumn{1}{c}{$8$} & \multicolumn{1}{c}{$13$} & \multicolumn{1}{c}{$21$} & \multicolumn{1}{c}{$34$} & \multicolumn{1}{c}{$55$} \\ 
 & Lucas & $L_r$ & \multicolumn{1}{c}{$7$} & \multicolumn{1}{c}{$-4$} & \multicolumn{1}{c}{$3$} & \multicolumn{1}{c}{$-1$} & \multicolumn{1}{c}{$2$} & \multicolumn{1}{c}{$1$} & \multicolumn{1}{c}{$3$} & \multicolumn{1}{c}{$4$} & \multicolumn{1}{c}{$7$} & \multicolumn{1}{c}{$11$} & \multicolumn{1}{c}{$18$} & \multicolumn{1}{c}{$29$} & \multicolumn{1}{c}{$47$} & \multicolumn{1}{c}{$76$} & \multicolumn{1}{c}{$123$} \\ 
\hline
$3$ & Tribonacci & $T_r$ & \multicolumn{1}{c}{$0$} & \multicolumn{1}{c}{$-1$} & \multicolumn{1}{c}{$1$} & \multicolumn{1}{c}{$0$} & \multicolumn{1}{c}{$0$} & \multicolumn{1}{c}{$1$} & \multicolumn{1}{c}{$1$} & \multicolumn{1}{c}{$2$} & \multicolumn{1}{c}{$4$} & \multicolumn{1}{c}{$7$} & \multicolumn{1}{c}{$13$} & \multicolumn{1}{c}{$24$} & \multicolumn{1}{c}{$44$} & \multicolumn{1}{c}{$81$} & \multicolumn{1}{c}{$149$} \\ 
 & Trib-Lucas & $K_r$ & \multicolumn{1}{c}{$-5$} & \multicolumn{1}{c}{$5$} & \multicolumn{1}{c}{$-1$} & \multicolumn{1}{c}{$-1$} & \multicolumn{1}{c}{$3$} & \multicolumn{1}{c}{$1$} & \multicolumn{1}{c}{$3$} & \multicolumn{1}{c}{$7$} & \multicolumn{1}{c}{$11$} & \multicolumn{1}{c}{$21$} & \multicolumn{1}{c}{$39$} & \multicolumn{1}{c}{$71$} & \multicolumn{1}{c}{$131$} & \multicolumn{1}{c}{$241$} & \multicolumn{1}{c}{$443$} \\ 
\hline
$4$ & Tetranacci & $M_r$ & \multicolumn{1}{c}{$-1$} & \multicolumn{1}{c}{$1$} & \multicolumn{1}{c}{$0$} & \multicolumn{1}{c}{$0$} & \multicolumn{1}{c}{$0$} & \multicolumn{1}{c}{$1$} & \multicolumn{1}{c}{$1$} & \multicolumn{1}{c}{$2$} & \multicolumn{1}{c}{$4$} & \multicolumn{1}{c}{$8$} & \multicolumn{1}{c}{$15$} & \multicolumn{1}{c}{$29$} & \multicolumn{1}{c}{$56$} & \multicolumn{1}{c}{$108$} & \multicolumn{1}{c}{$208$} \\ 
 & Tetra-Lucas & $R_r$ & \multicolumn{1}{c}{$7$} & \multicolumn{1}{c}{$-1$} & \multicolumn{1}{c}{$-1$} & \multicolumn{1}{c}{$-1$} & \multicolumn{1}{c}{$4$} & \multicolumn{1}{c}{$1$} & \multicolumn{1}{c}{$3$} & \multicolumn{1}{c}{$7$} & \multicolumn{1}{c}{$15$} & \multicolumn{1}{c}{$26$} & \multicolumn{1}{c}{$51$} & \multicolumn{1}{c}{$99$} & \multicolumn{1}{c}{$191$} & \multicolumn{1}{c}{$367$} & \multicolumn{1}{c}{$708$} \\ 
\hline
$5$ & Pentanacci & $P_r$ & \multicolumn{1}{c}{$1$} & \multicolumn{1}{c}{$0$} & \multicolumn{1}{c}{$0$} & \multicolumn{1}{c}{$0$} & \multicolumn{1}{c}{$0$} & \multicolumn{1}{c}{$1$} & \multicolumn{1}{c}{$1$} & \multicolumn{1}{c}{$2$} & \multicolumn{1}{c}{$4$} & \multicolumn{1}{c}{$8$} & \multicolumn{1}{c}{$16$} & \multicolumn{1}{c}{$31$} & \multicolumn{1}{c}{$61$} & \multicolumn{1}{c}{$120$} & \multicolumn{1}{c}{$236$} \\ 
 & Penta-Lucas & $Q_r$ & \multicolumn{1}{c}{$-1$} & \multicolumn{1}{c}{$-1$} & \multicolumn{1}{c}{$-1$} & \multicolumn{1}{c}{$-1$} & \multicolumn{1}{c}{$5$} & \multicolumn{1}{c}{$1$} & \multicolumn{1}{c}{$3$} & \multicolumn{1}{c}{$7$} & \multicolumn{1}{c}{$15$} & \multicolumn{1}{c}{$31$} & \multicolumn{1}{c}{$57$} & \multicolumn{1}{c}{$113$} & \multicolumn{1}{c}{$223$} & \multicolumn{1}{c}{$439$} & \multicolumn{1}{c}{$863$} \\ 
\end{tabular}
\caption{The first few sequences of the $n$-step Fibonacci numbers and $n$-step Lucas numbers.}
\end{table}
\clearpage
Our aim in writing this paper is to discover various properties of the generalized $n$-step Fibonacci numbers, $W_r$. Specifically we will develop recurrence relations, ordinary, binomial and double binomial summation identities, partial sums and generating functions.
\section{Recurrence relations}
\begin{theorem}\label{thm.xa5drh1}
The following identity holds, where $r$ and $s$ are integers:
\[
W_{r + s}=\sum_{i = 1}^n {\left( {\sum_{j = 0}^{n - i} {U_{s - j + 1} } } \right)W_{r - i} }\,.
\]
\end{theorem}
In particular, we have
\begin{equation}
U_{r + s}  = \sum_{i = 1}^n {\left( {\sum_{j = 0}^{n - i} {U_{s - j + 1} } } \right)U_{r - i} }
\end{equation}
and
\begin{equation}
V_{r + s}  = \sum_{i = 1}^n {\left( {\sum_{j = 0}^{n - i} {U_{s - j + 1} } } \right)V_{r - i} }\,.
\end{equation}
\begin{proof}
We will keep $r$ fixed and use induction on $s$. 

\medskip

The identity is true for $s=0$ because
\[
\sum_{j = 0}^{n - i} {U_{ - j + 1} }  = \sum_{j = i}^n {U_{ - n + j + 1} }  = \sum_{j = i}^{n - 1} {U_{ - n + j + 1} }  + 1 = 1\,,
\]
for $1\le i\le n$, by virtue of the initial terms \eqref{eq.c8583mo}.

\medskip

Assume that the identity is true for some integer $s=k\in\Z^+$. Let
\begin{equation}\label{eq.x1cggbs}
P_k: \left( W_{r + k}  =\sum_{i = 1}^n {\left( {\sum_{j = 0}^{n - i} {U_{k - j + 1} } } \right)W_{r - i} }\right)\,.
\end{equation}
We wish to prove that
\begin{equation}\label{eq.katdl81}
P_{k+1}:  \left(W_{r + k + 1}  = \sum_{i = 1}^n {\left( {\sum_{j = 0}^{n - i} {U_{k + 1 - j + 1} } } \right)W_{r - i} }\right)
\end{equation}
and
\begin{equation}\label{eq.an717w3}
P_{k-1}:  \left(W_{r + k - 1}  = \sum_{i = 1}^n {\left( {\sum_{j = 0}^{n - i} {U_{k - 1 - j + 1} } } \right)W_{r - i} }\right) 
\end{equation}
are true whenever $P_k$ holds.

\medskip

By the identity~\eqref{eq.ah0hczd} and the induction hypothesis $P_k$~(identity~\eqref{eq.x1cggbs}) we have
\begin{equation}\label{eq.ycv4kvc}
\begin{split}
W_{r + k + 1}  = \sum_{\lambda  = 1}^n {W_{r + k + 1 - \lambda } }  &= \sum_{\lambda  = 1}^n {\left\{ {\sum_{i = 1}^n {\left( {\sum_{j = 0}^{n - i} {U_{k + 1 - \lambda  - j + 1} } } \right)W_{r - i} } } \right\}}\\
&=\sum_{i = 1}^n {\sum_{j = 0}^{n - i} {\left( {\sum_{\lambda  = 1}^n {U_{k + 1 - \lambda  - j + 1} } } \right)} W_{r - i} }\,.
\end{split}
\end{equation}
By the recurrence relation~\eqref{eq.qnse2iy} we have
\begin{equation}\label{eq.weu8gnr}
\sum_{\lambda  = 1}^n {U_{k + 1 - \lambda  - j + 1} }  = U_{k + 1 - j + 1}\,.
\end{equation}
Using~\eqref{eq.weu8gnr} in~\eqref{eq.ycv4kvc} yields~\eqref{eq.katdl81} and therefore $P_k\Rightarrow P_{k+1}$. Following the same procedure, it is readily established that $P_k\Rightarrow P_{k-1}$. 
\end{proof}
We remark that Gabai \cite[Theorem 6]{gabai70} earlier proved the equivalent of Theorem \ref{thm.xa5drh1}. His proof, however, placed a restriction on the integers $r$ and $s$, consistent with his definition of the generalized $n$-step numbers. 
\begin{cor}\label{cor.p3a64yf}
The following identity holds, where $r$ and $s$ are integers:
\[
W_{r + s}=\sum_{i = 1}^n {\left( {\sum_{j = 0}^{n - i} {W_{s - j + 1} } } \right)U_{r - i} }\,.
\]

\end{cor}
In particular,
\begin{equation}\label{eq.tmq5485}
V_{r + s}=\sum_{i = 1}^n {\left( {\sum_{j = 0}^{n - i} {V_{s - j + 1} } } \right)U_{r - i} }\,.
\end{equation}
\begin{proof}
We require the following summation identities:
\begin{equation}\label{eq.dy52utj}
\sum_{j = a}^{k - a} {f_j }  = \sum_{j = a}^{k - a} {f_{k - j} }
\end{equation}
and
\begin{equation}\label{eq.yp43new}
\sum_{i = a}^n {\sum_{j = 0}^{n - i} {A_{i,i + j} } }  = \sum_{i = a}^n {\sum_{j = a}^i {A_{j,i} } }\,.
\end{equation}
Now,
\begin{equation}\label{eq.i6q706z}
W_{r + s}  = \sum_{i = 1}^n {\sum_{j = 0}^{n - i} {U_{s - j + 1} W_{r - i} } }  = \sum_{i = 1}^n {\sum_{j = 0}^{n - i} {W_{r - i} U_{s - j + 1} } }  = \sum_{i = 1}^n {\sum_{j = 0}^{n - i} {W_{r - i} U_{s - n + i + j + 1} } }\,,
\end{equation}
by application of identity \eqref{eq.dy52utj} to the $j$ summation.
Using identity \eqref{eq.yp43new} to re-write the sum in \eqref{eq.i6q706z} gives
\begin{equation}
W_{r + s}  = \sum_{i = 1}^n {\sum_{j = 1}^i {W_{r - j} U_{s - n + i + 1} } }\,,
\end{equation}
in which the application of identity \eqref{eq.dy52utj} to the $i$ summation gives
\begin{equation}\label{eq.vitr56w}
W_{r + s}  = \sum_{i = 1}^n {\sum_{j = 1}^{n + 1 - i} {W_{r - j} U_{s - n + n + 1 - i + 1} } }  = \sum_{i = 1}^n {\sum_{j = 1}^{n + 1 - i} {W_{r - j} U_{s - i + 2} } }  = \sum_{i = 1}^n {\sum_{j = 0}^{n - i} {W_{r - j - 1} U_{s - i + 2} } }\,.
\end{equation}
Finally, setting $r=s+2$ and $s=r-2$ in \eqref{eq.vitr56w} gives the identity of Corollary \ref{cor.p3a64yf}.
\end{proof}
We now give explicit examples of the identities of Theorem~\ref{thm.xa5drh1} and Corollary \ref{cor.p3a64yf} for low~$n$ $n$-step generalized Fibonacci numbers.
\subsection{Recurrence relations for the generalized Fibonacci numbers}
With $n=2$ in the identity of Theorem~\ref{thm.xa5drh1}, we have
\begin{equation}\label{eq.o0bcott}
\mathcal{F}_{r + s}  = F_{s + 2} \mathcal{F}_{r - 1}  + F_{s + 1} \mathcal{F}_{r - 2}\,,
\end{equation}
which is a variant of Formula~(8) of Vajda~\cite{vajda}, with particular instances
\begin{equation}
F_{r + s}  = F_{s + 2} F_{r - 1}  + F_{s + 1} F_{r - 2}
\end{equation}
and
\begin{equation}
L_{r + s}  = F_{s + 2} L_{r - 1}  + F_{s + 1} L_{r - 2}\,.
\end{equation}
\subsection{Recurrence relations for the generalized Tribonacci numbers}
Choosing $n=3$ in the identity of Theorem~\ref{thm.xa5drh1} gives
\begin{equation}\label{eq.oro9z3n}
\mathcal{T}_{r + s}  = T_{s + 2} \mathcal{T}_{r - 1}  + (T_{s + 1}  + T_s )\mathcal{T}_{r - 2}  + T_{s + 1} \mathcal{T}_{r - 3}\,,
\end{equation}
with the particular cases
\begin{equation}\label{eq.pc034en}
T_{r + s}  = T_{s + 2} T_{r - 1}  + (T_{s + 1}  + T_s )T_{r - 2}  + T_{s + 1} T_{r - 3}
\end{equation}
and
\begin{equation}\label{eq.uxvinyh}
K_{r + s}  = T_{s + 2} K_{r - 1}  + (T_{s + 1}  + T_s )K_{r - 2}  + T_{s + 1} K_{r - 3}\,.
\end{equation}
The identity~\eqref{eq.pc034en} was also proved by Feng~\cite{feng11} and by Shah~\cite{shah11}.

\medskip

Since $T_{-17}=0$, $T_{-18}=-103$ and $T_{-19}=159$, setting $s=-19$ in identity~\eqref{eq.oro9z3n} produces another three-term recurrence for the generalized Tribonacci numbers, namely
\begin{equation}\label{eq.mkpi1l0}
\mathcal{T}_{r-19}=56\mathcal{T}_{r-2}-103\mathcal{T}_{r-3}\,,
\end{equation}
in addition to the relation 
\begin{equation}
\mathcal{T}_{r}=2\mathcal{T}_{r-1}-\mathcal{T}_{r-4}\,,
\end{equation}
obtained at $n=3$ in identity~\eqref{eq.b5q12ro}.

\medskip

Choosing $n=3$ in the identity of Corollary \eqref{cor.p3a64yf} with $W=K$, $U=T$ gives
\begin{equation}\label{eq.j254zuq}
K_{r + s}  = K_{s + 2} T_{r - 1}  + (K_{s + 1}  + K_s )T_{r - 2}  + K_{s + 1} T_{r - 3}\,.
\end{equation}
Setting $s=-4$ in identity~\eqref{eq.j254zuq} gives a three-term identity connecting the Tribonacci-Lucas numbers and the Tribonacci numbers:
\begin{equation}\label{eq.ibzzex9}
K_{r-4}=-T_{r-1}+5T_{r-3}\,,
\end{equation}
since $K_{-3}=-K_{-4}=5$.
\subsection{Recurrence relations for the generalized Tetranacci numbers}
The choice $n=4$ in the identity of Theorem~\ref{thm.xa5drh1} gives
\begin{equation}
\begin{split}
\mathcal{M}_{r + s}  &= M_{s + 2} \mathcal{M}_{r - 1}  + (M_{s + 1}  + M_s  + M_{s - 1} )\mathcal{M}_{r - 2}\\
&\qquad + (M_{s + 1}  + M_s )\mathcal{M}_{r - 3}  + M_{s + 1} \mathcal{M}_{r - 4}\,,
\end{split}
\end{equation}
with the special cases
\begin{equation}
\begin{split}
M_{r + s}  &= M_{s + 2} M_{r - 1}  + (M_{s + 1}  + M_s  + M_{s - 1} )M_{r - 2}\\
&\qquad + (M_{s + 1}  + M_s )M_{r - 3}  + M_{s + 1} M_{r - 4}
\end{split}
\end{equation}
and
\begin{equation}
\begin{split}
R_{r + s}  &= M_{s + 2} R_{r - 1}  + (M_{s + 1}  + M_s  + M_{s - 1} )R_{r - 2}\\
&\qquad + (M_{s + 1}  + M_s )R_{r - 3}  + M_{s + 1} R_{r - 4}\,.
\end{split}
\end{equation}
Choosing $n=4$ in the identity of Corollary \eqref{cor.p3a64yf} with $W=R$, $U=M$ gives
\begin{equation}\label{eq.vny24q2}
\begin{split}
R_{r + s}  &= R_{s + 2} M_{r - 1}  + (R_{s + 1}  + R_s  + R_{s - 1} )M_{r - 2}\\
&\qquad + (R_{s + 1}  + R_s )M_{r - 3}  + R_{s + 1} M_{r - 4}\,.
\end{split}
\end{equation}
Setting $s=-9$, $s=-5$ and $s=-4$, respectively, in \eqref{eq.vny24q2}, yields, in each case, a four-term relation expressing a Tetranacci-Lucas number in terms of Tetranacci numbers:
\begin{equation}
R_{r - 9}  =  - M_{r - 1}  - 4M_{r - 3}  + 15M_{r - 4}\,,
\end{equation}
\begin{equation}
R_{r - 5}  =  - M_{r - 1}  + M_{r - 3}  + 7M_{r - 4}\,,
\end{equation}
\begin{equation}
R_{r - 4}  =  - M_{r - 1}  + 6M_{r - 3}  - M_{r - 4}\,.
\end{equation}
\section{Summation identities}
\begin{lemma}\label{lem.ewdmeif}
Let
\[
Z_r  = \sum_{j = 1}^{\left\lceil {n/2} \right\rceil } {W_{r - 2j + 1} }=
\begin{cases}
W_{r - 1}  + W_{r - 3}  + W_{r - 5}  +  \cdots  + W_{r - n+1}, & \text{if $n$ is even};\\
W_{r - 1}  + W_{r - 3}  + W_{r - 5}  +  \cdots  + W_{r - n}, & \text{if $n$ is odd},
\end{cases}
\]
where $\lceil q\rceil$ is the smallest integer greater than $q$.
Then
\begin{equation}
\begin{split}
Z_r+Z_{r-1}&=W_r+(n\bmod 2)W_{r-n-1}\\
&=
\begin{cases}
W_r,& \text{if $n$ is even};\\
2W_{r-1}, & \text{if $n$ is odd}.
\end{cases}
\end{split}
\end{equation}

\end{lemma}
\begin{lemma}[{\cite[Lemma 1]{adegoke18}}]\label{lem.u4bqbkc}
Let $\{X_r\}$ and $\{Y_r\}$ be any two sequences such that $X_r$ and $Y_r$, $r\in\Z$, are connected by a three-term recurrence relation $X_r=f_1X_{r-a}+f_2Y_{r-b}$, where $f_1$ and $f_2$ are arbitrary non-vanishing complex functions, not dependent on $r$, and $a$ and $b$ are integers. Then,
\[
f_2\sum_{j = 0}^k {\frac{{Y_{r - ka  - b  + a j} }}{{f_1{}^j }}}  = \frac{{X_r }}{{f_1{}^k }} - f_1X_{r - (k + 1)a }\,, 
\]
for $k$ a non-negative integer.
\end{lemma}
The next theorem follows directly from Lemma \ref{lem.ewdmeif} and Lemma \ref{lem.u4bqbkc}.
\begin{theorem}
The following identity holds, where $r$ and $k$ are integers:
\[
\begin{split}
&\sum_{j = 0}^k {( - 1)^j W_{r - k + j} }  + n\bmod 2\sum_{j = 0}^k {( - 1)^j W_{r - k - n - 1 + j} }\\
&\qquad= ( - 1)^k \sum_{j = 1}^{\left\lceil {n/2} \right\rceil } {W_{r - 2j + 1} }  + \sum_{j = 1}^{\left\lceil {n/2} \right\rceil } {W_{r - 2j - k} }\,.
\end{split}
\]

\end{theorem}
In particular,
\begin{equation}
\begin{split}
&\sum_{j = 0}^k {( - 1)^j U_{r - k + j} }  + n\bmod 2\sum_{j = 0}^k {( - 1)^j U_{r - k - 1 + j} }\\
&\qquad= ( - 1)^k \sum_{j = 1}^{\left\lceil {n/2} \right\rceil } {U_{r - 2j + 1} }  + \sum_{j = 1}^{\left\lceil {n/2} \right\rceil } {U_{r - 2j - k} }
\end{split}
\end{equation}
and
\begin{equation}
\begin{split}
&\sum_{j = 0}^k {( - 1)^j V_{r - k + j} }  + n\bmod 2\sum_{j = 0}^k {( - 1)^j V_{r - k - 1 + j} }\\
&\qquad= ( - 1)^k \sum_{j = 1}^{\left\lceil {n/2} \right\rceil } {V_{r - 2j + 1} }  + \sum_{j = 1}^{\left\lceil {n/2} \right\rceil } {V_{r - 2j - k} }\,.
\end{split}
\end{equation}
Thus, if $n$ is even, we have
\begin{equation}
\sum_{j = 0}^k {( - 1)^j W_{r - k + j} }  = ( - 1)^k \sum_{j = 1}^{n/2} {W_{r - 2j + 1} }  + \sum_{j = 1}^{n/2} {W_{r - 2j - k} }\,,
\end{equation}
while if $n$ is odd, we have
\begin{equation}
2\sum_{j = 0}^k {( - 1)^j W_{r - k + j - 1} }  = ( - 1)^k \sum_{j = 1}^{(n + 1)/2} {W_{r - 2j + 1} }  + \sum_{j = 1}^{(n + 1)/2} {W_{r - 2j - k} }\,.
\end{equation}
We give explicit examples with small $n$ values.
\begin{equation}
\sum_{j = 0}^k {( - 1)^j \mathcal{F}_{r - k + j} }  = ( - 1)^k \mathcal{F}_{r - 1}  + \mathcal{F}_{r - k - 2}\,,
\end{equation}
\begin{equation}
2\sum_{j = 0}^k {( - 1)^j \mathcal{T}_{r - k - 1 + j} }  = ( - 1)^k (\mathcal{T}_{r - 1}  + \mathcal{T}_{r - 3} ) + \mathcal{T}_{r - k - 2}  + \mathcal{T}_{r - k - 4}\,,
\end{equation}
\begin{equation}
\sum_{j = 0}^k {( - 1)^j \mathcal{M}_{r - k + j} }  = ( - 1)^k (\mathcal{M}_{r - 1}  + \mathcal{M}_{r - 3} ) + \mathcal{M}_{r - k - 2}  + \mathcal{M}_{r - k - 4}\,,
\end{equation}
\begin{equation}
2\sum_{j = 0}^k {( - 1)^j \mathcal{P}_{r - k - 1 + j} }  = ( - 1)^k (\mathcal{P}_{r - 1}  + \mathcal{P}_{r - 3}  + \mathcal{P}_{r - 5} ) + \mathcal{P}_{r - k - 2}  + \mathcal{P}_{r - k - 4}  + \mathcal{P}_{r - k - 6}\,.
\end{equation}
In particular,
\begin{equation}
\sum_{j = 0}^k {( - 1)^j \mathcal{F}_{j} }  = ( - 1)^k \mathcal{F}_{k - 1}  + \mathcal{F}_{ - 2}\,,
\end{equation}
\begin{equation}
2\sum_{j = 0}^k {( - 1)^j \mathcal{T}_{j} }  = ( - 1)^k (\mathcal{T}_{k}  + \mathcal{T}_{k - 2} ) + \mathcal{T}_{- 1}  + \mathcal{T}_{- 3}\,,
\end{equation}
\begin{equation}
\sum_{j = 0}^k {( - 1)^j \mathcal{M}_{j} }  = ( - 1)^k (\mathcal{M}_{k - 1}  + \mathcal{M}_{k - 3} ) + \mathcal{M}_{ - 2}  + \mathcal{M}_{ - 4}\,,
\end{equation}
\begin{equation}
2\sum_{j = 0}^k {( - 1)^j \mathcal{P}_j }  = ( - 1)^k (\mathcal{P}_k  + \mathcal{P}_{k - 2}  + \mathcal{P}_{r - 4} ) + \mathcal{P}_{ - 1}  + \mathcal{P}_{ - 3}  + \mathcal{P}_{ - 5}\,.
\end{equation}
\begin{lemma}[{\cite[Lemma 2]{adegoke18}}]\label{lem.s9jfs7n}
Let $\{X_r\}$ be any arbitrary sequence, where $X_r$, $r\in\Z$, satisfies a three-term recurrence relation $X_r=f_1X_{r-a}+f_2X_{r-b}$, where $f_1$ and $f_2$ are arbitrary non-vanishing complex functions, not dependent on $r$, and $a$ and $b$ are integers. Then, the following identities hold for integer $k$:
\begin{equation}\label{eq.mxyb9zk}
f_2\sum_{j = 0}^k {\frac{{X_{r - ka  - b  + a j} }}{{f_1^j }}}  = \frac{{X_r }}{{f_1^k }} - f_1X_{r - (k + 1)a }\,,
\end{equation}
\begin{equation}\label{eq.cgldajj}
f_1\sum_{j = 0}^k {\frac{{X_{r - kb  - a  + b j} }}{{f_2^j }}}  = \frac{{X_r }}{{f_2^k }} - f_2X_{r - (k + 1)b } 
\end{equation}
and
\begin{equation}\label{eq.n2n4ec3}
\sum_{j = 0}^k { \frac{X_{r - (a - b)k + b + (a - b)j}}{(-f_1/f_2)^j} }  = \frac{f_2 X_r}{(-f_1/f_2)^k}  + f_1 X_{r - (k + 1)(a - b)}\,.
\end{equation}

\end{lemma}
The next theorem is a consequence of identity~\eqref{eq.b5q12ro} and Lemma \ref{lem.s9jfs7n}.
\begin{theorem}\label{thm.sp0to1x}
The following identities hold, where $r$ and $k$ are integers:
\begin{equation}\label{eq.t0s0wfn}
\sum_{j = 0}^k {2^{k - j} W_{r - k - n - 1 + j} }  = 2^{k + 1} W_{r - k - 1}  - W_r\,,
\end{equation}
\begin{equation}
2\sum_{j = 0}^k {( - 1)^j W_{r - nk - k - 1 + (n + 1)j} }  = ( - 1)^k W_r  + W_{r - (k + 1)(n + 1)}
\end{equation}
and
\begin{equation}\label{eq.hytvkqp}
\sum_{j = 0}^k {2^j W_{r - nk + 1 + nj} }  = 2^{k + 1} W_r  - W_{r - (k + 1)n}\,.
\end{equation}

\end{theorem}
In particular,
\begin{equation}\label{eq.noe5wpc}
\sum_{j = 0}^k {2^{k - j} W_j }  = 2^{k + 1} W_n  - W_{k + n + 1}\,, 
\end{equation}
\begin{equation}
2\sum_{j = 0}^k {( - 1)^j W_{(n + 1)j} }  = ( - 1)^k W_{k(n + 1) + 1}  + 2W_0  - W_1 
\end{equation}
and
\begin{equation}
\sum_{j = 0}^k {2^j W_{nj} }  = 2^{k + 1} W_{kn - 1}  - 4W_{n - 1}  + 2W_n  + W_0\,. 
\end{equation}
We now illustrate Theorem \ref{thm.sp0to1x} for small values of $n$.
\subsection{Summation identities involving the generalized Fibonacci numbers, ($n=2$) }
\begin{equation}\label{eq.r409blq}
\sum_{j = 0}^k {2^{k - j} \mathcal{F}_{r - k - 3 + j} }  = 2^{k + 1} \mathcal{F}_{r - k - 1}  - \mathcal{F}_r\,,
\end{equation}
\begin{equation}
2\sum_{j = 0}^k {( - 1)^j \mathcal{F}_{r - 2k - k - 1 + 3j} }  = ( - 1)^k \mathcal{F}_r  + \mathcal{F}_{r - 3(k + 1)}
\end{equation}
and
\begin{equation}
\sum_{j = 0}^k {2^j \mathcal{F}_{r - 2k + 1 + 2j} }  = 2^{k + 1} \mathcal{F}_r  - \mathcal{F}_{r - 2(k + 1)}\,.
\end{equation}
In particular,
\begin{equation}
\sum_{j = 0}^k {2^{k - j} \mathcal{F}_j }  = 2^{k + 1} \mathcal{F}_2  - \mathcal{F}_{k + 3}\,, 
\end{equation}
\begin{equation}
2\sum_{j = 0}^k {( - 1)^j \mathcal{F}_{3j} }  = ( - 1)^k \mathcal{F}_{3k + 1}  + 2\mathcal{F}_0  - \mathcal{F}_1 
\end{equation}
and
\begin{equation}\label{eq.lx11b5k}
\sum_{j = 0}^k {2^j \mathcal{F}_{2j} }  = 2^{k + 1} \mathcal{F}_{2k - 1}  - 4\mathcal{F}_{1}  + 2\mathcal{F}_2  + \mathcal{F}_0\,. 
\end{equation}

\subsection{Summation identities involving the generalized Tribonacci numbers, ($n=3$) }
\begin{equation}
\sum_{j = 0}^k {2^{k - j} \mathcal{T}_{r - k - 4 + j} }  = 2^{k + 1} \mathcal{T}_{r - k - 1}  - \mathcal{T}_r\,,
\end{equation}
\begin{equation}
2\sum_{j = 0}^k {( - 1)^j \mathcal{T}_{r - 4k - 1 + 4j} }  = ( - 1)^k \mathcal{T}_r  + \mathcal{T}_{r - 4(k + 1)}
\end{equation}
and
\begin{equation}
\sum_{j = 0}^k {2^j \mathcal{T}_{r - 3k + 1 + 3j} }  = 2^{k + 1} \mathcal{T}_r  - \mathcal{T}_{r - 3(k + 1)}\,.
\end{equation}
In particular,
\begin{equation}
\sum_{j = 0}^k {2^{k - j} \mathcal{T}_j }  = 2^{k + 1} \mathcal{T}_3  - \mathcal{T}_{k + 4}\,, 
\end{equation}
\begin{equation}
2\sum_{j = 0}^k {( - 1)^j \mathcal{T}_{4j} }  = ( - 1)^k \mathcal{T}_{4k + 1}  + 2\mathcal{T}_0  - \mathcal{T}_1 
\end{equation}
and
\begin{equation}
\sum_{j = 0}^k {2^j \mathcal{T}_{3j} }  = 2^{k + 1} \mathcal{T}_{3k - 1}  - 4\mathcal{T}_{2}  + 2\mathcal{T}_3  + \mathcal{T}_0\,. 
\end{equation}
\subsection{Summation identities involving the generalized Tetranacci numbers, ($n=4$) }
\begin{equation}
\sum_{j = 0}^k {2^{k - j} \mathcal{M}_{r - k - 5 + j} }  = 2^{k + 1} \mathcal{M}_{r - k - 1}  - \mathcal{M}_r\,,
\end{equation}
\begin{equation}
2\sum_{j = 0}^k {( - 1)^j \mathcal{M}_{r - 4k - k - 1 + 5j} }  = ( - 1)^k \mathcal{M}_r  + \mathcal{M}_{r - 5(k + 1)}
\end{equation}
and
\begin{equation}
\sum_{j = 0}^k {2^j \mathcal{M}_{r - 4k + 1 + 4j} }  = 2^{k + 1} \mathcal{M}_r  - \mathcal{M}_{r - 4(k + 1)}\,.
\end{equation}
In particular,
\begin{equation}
\sum_{j = 0}^k {2^{k - j} \mathcal{M}_j }  = 2^{k + 1} \mathcal{M}_4  - \mathcal{M}_{k + 5}\,, 
\end{equation}
\begin{equation}
2\sum_{j = 0}^k {( - 1)^j \mathcal{M}_{5j} }  = ( - 1)^k \mathcal{M}_{5k + 1}  + 2\mathcal{M}_0  - \mathcal{M}_1 
\end{equation}
and
\begin{equation}
\sum_{j = 0}^k {2^j \mathcal{M}_{4j} }  = 2^{k + 1} \mathcal{M}_{4k - 1}  - 4\mathcal{M}_{3}  + 2\mathcal{M}_4  + \mathcal{M}_0\,. 
\end{equation}
\subsection{Further summation identities involving the generalized Fibonacci numbers}
In addition to the summation identities \mbox{\eqref{eq.r409blq} -- \eqref{eq.lx11b5k}}, we also have the results stated in the next theorem, on account of identity \eqref{eq.o0bcott} and \mbox{Lemma \ref{lem.s9jfs7n}}.
\begin{theorem}
The following identities hold, where $r$ and $s$ are integers:
\begin{equation}
F_s \sum_{j = 0}^k {F_{s + 1}^{k - j} \mathcal{F}_{r - 1 + sj} }  = \mathcal{F}_{r + s(k + 1)}  - F_{s + 1}^{k + 1} \mathcal{F}_r\,, 
\end{equation}
\begin{equation}
\sum_{j = 0}^k {( - 1)^j F_s^{k - j} F_{s + 1}^j \mathcal{F}_{r - k + s + j} }  = ( - 1)^k F_{s + 1}^{k + 1} \mathcal{F}_r  + F_s^{k + 1} \mathcal{F}_{r - k - 1} 
\end{equation}
and
\begin{equation}
F_s \sum_{j = 0}^k {F_{s - 1}^{k - j} \mathcal{F}_{r - sk - s + 1 + sj} }  = \mathcal{F}_r  - F_{s - 1}^{k + 1} \mathcal{F}_{r - (k + 1)s}\,.
\end{equation}
In particular,
\begin{equation}
F_s \sum_{j = 0}^k {F_{s + 1}^{k - j} \mathcal{F}_{sj} }  = \mathcal{F}_{sk + s + 1}  - F_{s + 1}^{k + 1} \mathcal{F}_1 \,,
\end{equation}
\begin{equation}
\sum_{j = 0}^k {( - 1)^j F_s^{k - j} F_{s + 1}^j \mathcal{F}_j }  = ( - 1)^k F_{s + 1}^{k + 1} \mathcal{F}_{k - s}  + F_s^{k + 1} \mathcal{F}_{ - s - 1} 
\end{equation}
and
\begin{equation}
F_s \sum_{j = 0}^k {F_{s - 1}^{k - j} \mathcal{F}_{sj} }  = \mathcal{F}_{sk + s - 1}  - F_{s - 1}^{k + 1} \mathcal{F}_{ - 1} \,.
\end{equation}

\end{theorem}
When identity \eqref{eq.o0bcott} is written as
\begin{equation}
\mathcal{F}_{s - 1} F_r  =  - \mathcal{F}_s F_{r + 1}  + \mathcal{F}_{r + s}
\end{equation}
and the identifications $X=F$ and $Y=\mathcal{F}$ are made in Lemma \ref{lem.u4bqbkc} we have the result stated in the next theorem.
\begin{theorem}
The following identity holds where $r$, $s$ and $k$ are integers:
\begin{equation}
\sum_{j = 0}^k {( - 1)^j \mathcal{F}_{s - 1}^{k - j} \mathcal{F}_s^j \mathcal{F}_{r + s + j} }  = F_r \mathcal{F}_{s - 1}^{k + 1}  - ( - 1)^{k + 1} F_{r + k + 1} \mathcal{F}_s^{k + 1}\,.
\end{equation}

\end{theorem}
In particular,
\begin{equation}
\sum_{j = 0}^k {( - 1)^j \mathcal{F}_{s - 1}^{k - j} \mathcal{F}_s^j \mathcal{F}_j }  = ( - 1)^{s - 1} F_s \mathcal{F}_{s - 1}^{k + 1}  - ( - 1)^{k - 1} F_{k - s + 1} \mathcal{F}_s^{k + 1}\,.
\end{equation}
\subsection{Further summation identities involving the generalized Tribonacci numbers}
The next theorem, expressing a summation involving Tribonacci-Lucas numbers in terms of Tribonacci numbers, follows from identity~\eqref{eq.ibzzex9} and Lemma \ref{lem.u4bqbkc}.
\begin{thm}
The following identity holds, where $r$ and $k$ are integers:
\[
\sum_{j = 0}^k {5^{k - j} K_{r - 2k - 3 + 2j} }  = T_r  - 5^{k + 1} T_{r - 2k - 2}\,. 
\]

\end{thm}
In particular,
\begin{equation}
\sum_{j = 0}^k {5^{k - j} K_{2j} }  = T_{2k + 3}  - 5^{k + 1}\,.
\end{equation}
Further summation identities are obtained from identity \eqref{eq.mkpi1l0} and Lemma \ref{lem.s9jfs7n}. These are presented in the next theorem.
\begin{theorem}
The following identities hold, where $r$ and $k$ are integers:
\begin{equation}
103\sum_{j = 0}^k {56^j \mathcal{T}_{r + 16 + 17j} }  = 56^{k + 1} \mathcal{T}_{r + 17k + 17}  - \mathcal{T}_r\,,
\end{equation}
\begin{equation}
56\sum_{j = 0}^k {( - 1)^j 103^j \mathcal{T}_{r + 17 + 16j} }  = \mathcal{T}_r  - ( - 103)^{k + 1} \mathcal{T}_{r + 16k + 16}
\end{equation}
and
\begin{equation}
\sum_{j = 0}^k {103^{k - j} 56^j \mathcal{T}_{r - 16 + j} }  =  - 103^{k + 1} \mathcal{T}_r  + 56^{k + 1} \mathcal{T}_{r + k + 1}\,.
\end{equation}

\end{theorem}
In particular,
\begin{equation}
103\sum_{j = 0}^k {56^j \mathcal{T}_{17j} }  = 56^{k + 1} \mathcal{T}_{17k + 1}  - \mathcal{T}_{ - 16} \,,
\end{equation}
\begin{equation}
56\sum_{j = 0}^k {( - 1)^j 103^j \mathcal{T}_{16j} }  = \mathcal{T}_{ - 17}  - ( - 103)^{k + 1} \mathcal{T}_{16k - 1} 
\end{equation}
and
\begin{equation}
\sum_{j = 0}^k {103^{k - j} 56^j \mathcal{T}_j }  =  - 103^{k + 1} \mathcal{T}_{16}  + 56^{k + 1} \mathcal{T}_{k + 17}\,.
\end{equation}
\section{Binomial summation identities}
\begin{lemma}[{\cite[Lemma 3]{adegoke18}}]\label{lem.binomial}
Let $\{X_r\}$ be any arbitrary sequence. Let $X_r$, $r\in\Z$, satisfy a three-term recurrence relation $X_r=f_1X_{r-a}+f_2X_{r-b}$, where $f_1$ and $f_2$ are non-vanishing complex functions, not dependent on $r$, and $a$ and $b$ are integers. Then,
\begin{equation}\label{eq.fe496kc}
\sum_{j = 0}^k {\binom kj\left( {\frac{f_1}{f_2}} \right)^j X_{r - bk  + (b  - a )j} }  = \frac{{X_r }}{{f_2^k }}\,,
\end{equation}
\begin{equation}\label{eq.j7k6a8g}
\sum_{j = 0}^k {\binom kj\frac{{X_{r + (a - b)k + bj} }}{{( - f_2 )^j }}}  = \left( { - \frac{{f_1 }}{{f_2 }}} \right)^k X_r
\end{equation}
and
\begin{equation}\label{eq.fnwrzi3}
\sum_{j = 0}^k {\binom kj\frac{{X_{r + (b - a)k + a j} }}{{( - f_1 )^j }}}  = \left( { - \frac{f_2}{f_1}} \right)^k X_r\,,
\end{equation}
for $k$ a non-negative integer.

\end{lemma}
The next theorem is a consequence of identity~\eqref{eq.b5q12ro} and Lemma \ref{lem.binomial}.
\begin{thm}\label{thm.binomial}
The following identities hold, where $k$ is any non-negative integer and $r$ is any integer:
\begin{equation}\label{eq.uzehcqz}
\sum_{j = 0}^k {( - 1)^j \binom kj2^j W_{r - (n + 1)k + nj} }  = ( - 1)^k W_r\,,
\end{equation}
\begin{equation}
\sum_{j = 0}^k {\binom kjW_{r - nk + (n + 1)j} }  = 2^k W_r 
\end{equation}
and
\begin{equation}\label{eq.pn3n7iq}
\sum_{j = 0}^k {( - 1)^j \binom kj2^{k - j} W_{r + nk + j} }  = W_r\,. 
\end{equation}

\end{thm}
In particular,
\begin{equation}
\sum_{j = 0}^k {( - 1)^j \binom kj2^j W_{nj} }  = ( - 1)^k W_{(n+1)k}\,,
\end{equation}
\begin{equation}
\sum_{j = 0}^k {\binom kjW_{(n + 1)j} }  = 2^k W_{nk} 
\end{equation}
and
\begin{equation}\label{eq.j3u95dt}
\sum_{j = 0}^k {( - 1)^j \binom kj2^{k - j} W_{j} }  = W_{-nk}\,. 
\end{equation}
We remark that identity \eqref{eq.j3u95dt} proves Conjecture 2 (equation (15)) of Hisert \cite{hisert15}.
\subsection{Further binomial summation identities involving generalized Fibonacci numbers}
In addition to the summation identities obtained by setting $n=2$ in identities \mbox{\eqref{eq.uzehcqz} -- \eqref{eq.j3u95dt}} of Theorem \ref{thm.binomial}, we also have the results stated in the next theorem, on account of identity \eqref{eq.o0bcott} and \mbox{Lemma \ref{lem.binomial}}.
\begin{thm}
The following identities hold, where $k$ is any non-negative integer and $r$ and $s$ are any integers:
\begin{equation}
\sum_{j = 0}^k {( - 1)^j \binom kjF_{s-1}^{k - j} \mathcal{F}_{r - k + sj} }  = ( - 1)^k F_{s}^k \mathcal{F}_r \,,
\end{equation}
\begin{equation}
\sum_{j = 0}^k {\binom kjF_{s-1}^{k - j} F_{s}^j \mathcal{F}_{r - sk + j} }  = \mathcal{F}_r
\end{equation}
and
\begin{equation}
\sum_{j = 0}^k {( - 1)^{k - j} \binom kjF_{s + 1}^{k - j} \mathcal{F}_{r + k + sj} }  = F_s^k \mathcal{F}_r\,.
\end{equation}

\end{thm}
In particular,
\begin{equation}
\sum_{j = 0}^k {( - 1)^j \binom kjF_{s-1}^{k - j} \mathcal{F}_{sj} }  = ( - 1)^k F_{s}^k \mathcal{F}_k \,,
\end{equation}
\begin{equation}
\sum_{j = 0}^k {\binom kjF_{s-1}^{k - j} F_{s}^j \mathcal{F}_{j} }  = \mathcal{F}_{sk}
\end{equation}
and
\begin{equation}
\sum_{j = 0}^k {( - 1)^{k - j} \binom kjF_{s + 1}^{k - j} \mathcal{F}_{sj} }  = F_s^k \mathcal{F}_{-k}\,.
\end{equation}
\subsection{Further binomial summation identities involving generalized Tribonacci numbers}
In addition to the summation identities obtained by setting $n=3$ in identities \mbox{\eqref{eq.uzehcqz} -- \eqref{eq.j3u95dt}} of Theorem \ref{thm.binomial}, we also have the results stated in the next theorem, on account of identity \eqref{eq.mkpi1l0} and \mbox{Lemma \ref{lem.binomial}}.
\begin{theorem}
The following identities hold, where $k$ and $r$ are integers:
\begin{equation}
\sum_{j = 0}^k {( - 1)^{k - j} \binom kj103^{k - j} 56^j \mathcal{T}_{r + 16k + j} }  = \mathcal{T}_r\,, 
\end{equation}
\begin{equation}
\sum_{j = 0}^k {103^j \binom kj\mathcal{T}_{r - 17k + 16j} }  = 56^k \mathcal{T}_r 
\end{equation}
and
\begin{equation}
\sum_{j = 0}^k {( - 1)^j \binom kj56^j \mathcal{T}_{r - 16k + 17j} }  = ( - 103)^k \mathcal{T}_r\,. 
\end{equation}

\end{theorem}
In particular,
\begin{equation}
\sum_{j = 0}^k {( - 1)^{k - j} \binom kj103^{k - j} 56^j \mathcal{T}_j }  = \mathcal{T}_{ - 16k} \,,
\end{equation}
\begin{equation}
\sum_{j = 0}^k {103^j \binom kj\mathcal{T}_{16j} }  = 56^k \mathcal{T}_{17k} 
\end{equation}
and
\begin{equation}
\sum_{j = 0}^k {( - 1)^j \binom kj56^j \mathcal{T}_{17j} }  = ( - 103)^k \mathcal{T}_{16k}\,.
\end{equation}
\section{Double binomial summation identities}
\begin{lemma}[{\cite[Lemma 5]{adegoke18c}}]\label{lem.h2de9i7}
Let $\{X_r\}$ be any arbitrary sequence, $X_r$ satisfying a four-term recurrence relation $X_r=f_1X_{r-a}+f_2X_{r-b}+f_3X_{r-c}$, where $f_1$, $f_2$ and $f_3$ are arbitrary nonvanishing functions and $a$, $b$ and $c$ are integers. Then, the following identities hold:
\begin{equation}\label{eq.wgx2r2f}
\sum_{j = 0}^k {\sum_{s = 0}^j {\binom kj\binom js\left( {\frac{{f_2 }}{{f_3 }}} \right)^j\left( {\frac{{f_1 }}{{f_2 }}} \right)^s  X_{r - ck + (c - b)j + (b - a)s} } }  = \frac{{X_r }}{{f_3{}^k }}\,,
\end{equation}
\begin{equation}\label{eq.sm9bygb}
\sum_{j = 0}^k {\sum_{s = 0}^j {\binom kj\binom js\left( {\frac{{f_3 }}{{f_2 }}} \right)^j\left( {\frac{{f_1 }}{{f_3 }}} \right)^s  X_{r - bk + (b - c)j + (c - a)s} } }  = \frac{{X_r }}{{f_2{}^k }}\,,
\end{equation}
\begin{equation}
\sum_{j = 0}^k {\sum_{s = 0}^j {\binom kj\binom js\left( {\frac{{f_3 }}{{f_1 }}} \right)^j\left( {\frac{{f_2 }}{{f_3 }}} \right)^s X_{r - ak + (a - c)j + (c - b)s} } }  = \frac{{X_r }}{{f_1{}^k }}\,,
\end{equation}
\begin{equation}
\sum_{j = 0}^k {\sum_{s = 0}^j {\binom kj\binom js\left( {\frac{{f_2 }}{{f_3 }}} \right)^j\left( {-\frac{{1 }}{{f_2 }}} \right)^s X_{r - (c-a)k + (c - b)j + bs } } }  = \left(-\frac {f_1}{f_3}\right)^kX_r\,,
\end{equation}
\begin{equation}
\sum_{j = 0}^k {\sum_{s = 0}^j {\binom kj\binom js\left( {\frac{{f_1 }}{{f_3 }}} \right)^j\left( {-\frac{{1 }}{{f_1 }}} \right)^s X_{r - (c-b)k + (c - a)j + as } } }  = \left(-\frac {f_2}{f_3}\right)^kX_r\,,
\end{equation}
and
\begin{equation}\label{eq.o7540wl}
\sum_{j = 0}^k {\sum_{s = 0}^j {\binom kj\binom js\left( {\frac{{f_1 }}{{f_2 }}} \right)^j\left( {-\frac{{1 }}{{f_1 }}} \right)^s X_{r  - (b-c)k + (b - a)j + as} } }  = \left(-\frac {f_3}{f_2}\right)^kX_r\,.
\end{equation}

\end{lemma}
Evaluating identities \eqref{eq.t0s0wfn}--\eqref{eq.hytvkqp} at $k=1$ produces the following recurrence relations:
\begin{equation}\label{eq.jchf898}
W_r  = 4W_{r - 2}  - W_{r - n - 1}  - 2W_{r - n - 2}\,,
\end{equation}
\begin{equation}
W_r  = 2W_{r - 1}  - 2W_{r - n - 2}  + W_{r - 2n - 2}
\end{equation}
and
\begin{equation}
2W_r  = 4W_{r - 1}  - W_{r - n}  - W_{r - 2n - 1}\,.
\end{equation}
Evaluating identities \eqref{eq.uzehcqz}--\eqref{eq.pn3n7iq} at $k=2$ gives the following recurrence relations:
\begin{equation}
W_r  = 4W_{r - 2}  - 4W_{r - n - 2}  + W_{r - 2n - 2}\,
\end{equation}
\begin{equation}
W_r  = 4W_{r - 2}  - 2W_{r - n - 1}  - W_{r - 2n - 2}
\end{equation}
and
\begin{equation}\label{eq.ndejisc}
W_r  = 4W_{r - 1}  - 4W_{r - 2}  + W_{r - 2n - 2}\,.
\end{equation}
Each of identities \eqref{eq.jchf898}--\eqref{eq.ndejisc} has six double binomial summation identities associated with it.
In the next theorem we give the double binomial summation identities resulting from identity \eqref{eq.ndejisc}.
\begin{theorem}
The following identities hold for nonnegative integer $k$ and any integer $r$:
\begin{equation}
\sum_{j = 0}^k {\sum_{s = 0}^j {( - 1)^{j + s} \binom kj\binom js4^j W_{r - (2n + 2)k + 2nj + s} } }  = W_r\,, 
\end{equation}
\begin{equation}
\sum_{j = 0}^k {\sum_{s = 0}^j {( - 4)^{k - j} \binom kj\binom js4^s W_{r - 2k - 2nj + (2n + 1)s} } }  = W_r\,, 
\end{equation}
\begin{equation}
\sum_{j = 0}^k {\sum_{s = 0}^j {( - 1)^s \binom kj\binom js4^{k - j + s} W_{r - k - (2n + 1)j + 2ns} } }  = W_r\,, 
\end{equation}
\begin{equation}
\sum_{j = 0}^k {\sum_{s = 0}^j {( - 1)^{j - k} \binom kj\binom js4^{j - k - s} W_{r - (2n + 1)k + 2nj + 2s} } }  = W_r\,, 
\end{equation}
\begin{equation}
\sum_{j = 0}^k {\sum_{s = 0}^j {( - 1)^s \binom kj\binom js4^{j - k - s} W_{r - 2nk + j + s} } }  = W_r 
\end{equation}
and
\begin{equation}
\sum_{j = 0}^k {\sum_{s = 0}^j {( - 1)^{j + s} \binom kj\binom js4^{k - s} W_{r + 2nk + j + s} } }  = W_r\,.
\end{equation}

\end{theorem}
\section{Partial sums and generating function}
\begin{lemma}[{\cite[Lemma 2]{adegoke18c}}Partial sum of a $n$-term sequence]\label{lemma.qm8k37h}
Let $\{X_j\}$ be any arbitrary sequence, where $X_j$, $j\in\Z$, satisfies a $n$-term recurrence relation $X_j=f_1X_{j-c_1}+f_2X_{j-c_2}+\cdots+f_nX_{j-c_n}=\sum_{m=1}^n f_mX_{j-c_m}$, where $f_1$, $f_2$, $\ldots$, $f_n$ are arbitrary non-vanishing complex functions, not dependent on $j$, and $c_1$, $c_2$, $\ldots$, $c_n$ are fixed integers. Then, the following summation identity holds for arbitrary $x$ and non-negative integer $k$ :
\[
\sum_{j = 0}^k {x^j X_j }  = \frac{{\sum_{m = 1}^n {\left\{ {x^{c_m } f_m \left( {\sum_{j = 1}^{c_m } {x^{ - j} X_{ - j} }  - \sum_{j = k - c_m  + 1}^k {x^j X_j } } \right)} \right\}} }}{{1 - \sum_{m = 1}^n {x^{c_m } f_m } }}\,.
\]

\end{lemma}
We note that a special case of Lemma~\ref{lemma.qm8k37h} was proved in~\cite{zeitlin64}.

\medskip

The next theorem follows directly from Lemma \ref{lemma.qm8k37h} on account of identity \eqref{eq.b5q12ro}.
\begin{theorem}\label{thm.c9yvw0a}
The following identity holds for $k$ an integer and any $x$:
\[
(1 - 2x + x^{n + 1} )\sum_{j = 0}^k {x^j W_j }  = 2W_{ - 1}  - 2x^{k + 1} W_k  + x^{n + 1} \sum_{j = k - n}^k {x^j W_j }  - x^{n + 1} \sum_{j = 1}^{n + 1} {x^{ - j} W_{ - j} }\,. 
\]
\end{theorem}
We now work out the special cases of the identity of Theorem \ref{thm.c9yvw0a} for the $n$-step Fibonacci and $n$-step Lucas numbers.

\medskip

Now,
\begin{equation}
\begin{split}
\sum_{j = 1}^{n + 1} {x^{ - j} U_{ - j} }  &= x^{ - 1} U_{ - 1}  + x^{ - 2} U_{ - 2}  +  \cdots  + x^{ - n + 2} U_{ - n + 2}\\
&\qquad+ x^{ - n + 1} U_{ - n + 1}  + x^{ - n} U_{ - n}  + x^{ - n - 1} U_{ - n - 1}\,. 
\end{split}
\end{equation}
All except the last three terms on the right hand side of the above expression vanish on account of the initial terms as given in equation \eqref{eq.c8583mo}. Thus,
\begin{equation}\label{eq.qg1zyr3}
\begin{split}
\sum_{j = 1}^{n + 1} {x^{ - j} U_{ - j} }  &= x^{ - n + 1} U_{ - n + 1}  + x^{ - n} U_{ - n}  + x^{ - n - 1} U_{ - n - 1}\\
&\quad= x^{ - n + 1}  - x^{ - n}  + 2\delta _{n,2} x^{ - n - 1},\quad\mbox{by \eqref{eq.c8583mo} and \eqref{eq.fz5izyy}}\,.
\end{split}
\end{equation}
Using \eqref{eq.qg1zyr3} in the identity of Theorem \ref{thm.c9yvw0a} with $W=U$ we have
\begin{equation}\label{eq.j8vwfco}
(1 - 2x + x^{n + 1} )\sum_{j = 0}^k {x^j U_j }  = x - x^2  - 2x^{k + 1} U_k  + x^{n + 1} \sum_{j = k - n}^k {x^j U_j }\,.
\end{equation}
Next, we find
\begin{equation}
\begin{split}
\sum_{j = 1}^{n + 1} {x^{ - j} V_{ - j} }  &= x^{ - 1} V_{ - 1}  + x^{ - 2} V_{ - 2}  +  \cdots  + x^{ - n + 1} V_{ - n + 1}  + x^{ - n} V_{ - n}  + x^{ - n - 1} V_{ - n - 1} \\
&=  - (x^{ - 1}  + x^{ - 2}  +  \cdots  + x^{ - n + 1} ) + (2n - 1)x^{ - n}  - (n + 2)x^{ - n - 1}\,;
\end{split}
\end{equation}
so that,
\begin{equation}\label{eq.y43q21h}
\begin{split}
x^{n + 1} \sum_{j = 1}^{n + 1} {x^{ - j} V_{ - j} }  &=  - (x^n  + x^{n - 1}  +  \cdots  + x^3  + x^2 ) + (2n - 1)x - (n + 2)\\
&=  - \frac{x^{n + 1} -x^2}{x - 1} + (2n - 1)x - (n + 2)\,.
\end{split}
\end{equation}
Putting \eqref{eq.y43q21h} in the identity of Theorem \ref{thm.c9yvw0a} with $W=V$ we have
\begin{equation}\label{eq.buwdmoo}
\begin{split}
(1 - x)(1 - 2x + x^{n + 1} )\sum_{j = 0}^k {x^j V_j }  &= n - (3n - 1)x + 2nx^2  - x^{n + 1}\\
&\quad- (1 - x)x^{k + 1} 2V_k  + (1 - x)x^{n + 1} \sum_{j = k - n}^k {x^j V_j }\,.
\end{split}
\end{equation}
Note that the identity of Theorem \ref{thm.c9yvw0a} cannot be used directly to compute $\sum_{j=0}^k{W_j}$ because
\begin{equation}
\sum_{j = k - n}^k {W_j }  = \sum_{j = 0}^n {W_{j + k - n} }  = \sum_{j = 0}^n {W_{k - j} }  = W_k  + \sum_{j = 1}^n {W_{k - j} }  = 2W_k
\end{equation}
and
\begin{equation}
\sum_{j = 1}^{n + 1} {W_{ - j} }  = \sum_{j = 1}^n {W_{ - j} }  + W_{ - n - 1}  = W_0  + W_{ - n - 1}  = 2W_{ - 1}\,;
\end{equation}
so that both sides of the identity of Theorem \ref{thm.c9yvw0a} evaluates to zero at $x=1$. Nevertheless, the said sum can be evaluated if we divide both sides of the identity by $1-2x+x^{n+1}$ and then use L'Hospital's rule to take the limit at $x=1$, giving
\begin{equation}
(n - 1)\sum_{j = 0}^k {W_j }  = 2(n - k)W_k  - 2(n + 1)W_{ - 1}  + \sum_{j = k - n}^k {jW_j }  + \sum_{j = 1}^{n + 1} {jW_{ - j} }\,.
\end{equation}
Since
\begin{equation}
\begin{split}
\sum_{j = 1}^{n + 1} {jU_{ - j} }  &= \sum_{j = 1}^{n - 2} {jU_{ - j} }  + (n - 1)U_{ - n + 1}  + nU_{ - n}  + (n + 1)U_{ - n - 1}\\
&= 2(n + 1)\delta _{n,2}  - 1\,;
\end{split}
\end{equation}
and
\begin{equation}
\begin{split}
\sum_{j = 1}^{n + 1} {jV_{ - j} }  &= \sum_{j = 1}^{n - 1} {jV_{ - j} }  + nV_{ - n}  + (n + 1)V_{ - n - 1}\\
&=  - \sum_{j = 1}^{n - 1} j  + nV_{ - n}  + (n + 1)V_{ - n - 1} \\
&= \frac{{n^2 }}{2} - \frac{{7n}}{2} - 2\,,
\end{split}
\end{equation}
we obtain the following results for the sum of the first $k+1$ terms of the $n$-step Fibonacci numbers and the first $k+1$ terms of the $n$-step Lucas numbers:
\begin{equation}
(n - 1)\sum_{j = 0}^k {U_j }  = -1 + 2(n - k)U_k + \sum_{j = k - n}^k {jU_j }
\end{equation}
and
\begin{equation}
2(n - 1)\sum_{j = 0}^k {V_j }  = n(n - 3) + 4(n - k)V_k  + 2\sum_{j = k - n}^k {jV_j }\,.
\end{equation}
\begin{lemma}[{\cite[Lemma 3]{adegoke18c}}Generating function]\label{lemma.v1j9biq}
Under the conditions of Lemma~\ref{lemma.qm8k37h}, if additionally $x^kX_k$ vanishes in the limit as $k$ approaches infinity, then
\[
G_X  (x) = \sum_{j = 0}^\infty  {x^j X_j }  = \frac{{\sum_{m = 1}^n {\left( {x^{c_m } f_m \sum_{j = 1}^{c_m } {x^{ - j} X_{ - j} } } \right)} }}{{1 - \sum_{m = 1}^n {x^{c_m } f_m } }}\,,
\]
so that $G_X(x)$ is a generating function for the sequence $\{X_j\}$.
\end{lemma}
\begin{theorem}
The generalized $n$-step Fibonacci numbers have the following generating funtion:
\[
G_W(x;n)=\sum_{j = 0}^\infty  {x^j W_j }  = \frac{{2W_{ - 1}  - x^{n + 1} \sum_{j = 1}^{n + 1} {x^{ - j} W_{ - j} } }}{{1 - 2x + x^{n + 1} }}\,.
\]
\end{theorem}
In particular, from \eqref{eq.j8vwfco} and \eqref{eq.buwdmoo}, we see that the $n-$step Fibonacci and $n-$step Lucas numbers are generated, respectively, by
\begin{equation}
G_U(x;n)=\sum_{j = 0}^\infty  {x^j U_j }  = \frac{{x(1 - x)}}{{1 - 2x + x^{n + 1} }}
\end{equation}
and
\begin{equation}
G_V(x;n)=\sum_{j = 0}^\infty  {x^j V_j }  = \frac{{n - (3n - 1)x + 2nx^2  - x^{n + 1} }}{(1-x)(1 - 2x + x^{n + 1}) }\,.
\end{equation}

\hrule

\noindent 2010 {\it Mathematics Subject Classification}:
Primary 11B39; Secondary 11B37.

\noindent \emph{Keywords: }
$n$-step Fibonacci number, $n$-step Lucas number, Fibonacci number, Tribonacci number, Tetranacci number, Pentanacci number, summation identity, recurrence relation, generating function, partial sum.

\hrule

\end{document}